\documentclass{amsart} 
\usepackage{amsfonts,amsmath,amsthm,amssymb}
\usepackage{epsfig}
\usepackage{latexsym}
\begin{document}

\title{The modular isomorphism problem for finite $p$-groups
       with a cyclic subgroup of index $p^2$}
\author{Czes{\l}aw Bagi\'nski and Alexander Konovalov}
\date{}

\maketitle

\newtheorem{theorem}{Theorem}
\newtheorem{lemma}{Lemma}
\newcommand{\cl}{\mathop{\rm cl}\nolimits}
\hyphenation{non-abelian}
\hyphenation{un-publi-shed}
\hyphenation{diplom-arbeit}
\hyphenation{appli-cation}
\hyphenation{packa-ge}
\pagestyle{plain}

\begin{abstract}

Let $p$ be a prime number, $G$ be a finite $p$-group and $K$ be a
field of characteristic \nolinebreak $p$. The Modular Isomorphism
Problem (MIP) asks whether the group algebra $KG$ determines the
group $G$.
Dealing with MIP, we investigated a question whether the
nilpotency class of a finite $p$-group is determined by its
modular group algebra over the field of $p$ elements. We give a
positive answer to this question provided one of the following
conditions holds:
(i) $\exp G=p$;
(ii) $\cl(G)=2$;
(iii) $G'$ is cyclic;
(iv) $G$ is a group of maximal class and contains an abelian
subgroup of index $p$.

As a consequence, the positive solution of MIP for all $p$-groups
containing a cyclic subgroup of index $p^2$ was obtained.
\end{abstract}

\section{Introduction}
\label{intro}

Though the Modular Isomorphism Problem is known for more than 50
years, up to now it remains open. It was solved only for some
classes of $p$-groups, in particular:
\begin{itemize}
\item abelian $p$-groups (Deskins \cite{Deskins};
      alternate proof by Coleman \cite{Coleman});
\item $p$-groups of class 2 with elementary abelian commutator
      subgroup (Sandling, theorem 6.25 in \cite{Sandling-survey});
\item metacyclic $p$-groups (for $p>3$ by Bagi\'nski \cite{Baginski-meta};
      completed by Sandling \cite{Sandling-meta});
\item 2-groups of maximal class (Carlson \cite{Carlson};
      alternate proof by Bagi\'nski \cite{Baginski-maxclass2});
\item $p$-groups of maximal class, $p \ne 2$, when $|G| \le p^{p+1}$
      and $G$ contains an abelian maximal subgroup
      (Caranti and Bagi\'nski \cite{CJM1988});
\item elementary abelian-by-cyclic groups (Bagi\'nski \cite{Baginski-elab});
\item $p$-groups with the center of index $p^2$ (Drensky \cite{Drensky}),
\end{itemize}
where the results for abelian case, 2-groups of maximal class and
$p$-groups with the center of index $p^2$ are valid for arbitrary fields
of characteristic $p$.
Also it was solved for a number of groups of small
orders and a field of $p$ elements, in particular:
\begin{itemize}
\item groups of order not greater then $p^4$ (Passman \cite{Passman});
\item groups of order $2^5$ (Makasikis \cite{Makasikis} with remarks by
      Sandling \cite{Sandling-survey}; alternate proof by Michler,
      Newman and O'Brien \cite{Michler});
\item groups of order $p^5$ (Kovacs and Newman, due to Sandling's remark
      in \cite{Sandling-elab}; alternate proof by Salim and Sandling
      \cite{Salim-p5, Salim-maxclass});
\item groups of order $2^6$ (Wursthorn, using
      computer \cite{Wursthorn-dip, Wursthorn});
\item groups of order $2^7$ (Wursthorn, using
      computer \cite{Bleher}).
\end{itemize}

Besides this, a lot of MIP invariants (i.e. group properties which
are determined by the group algebra) are known, and they are very
useful for researches in MIP. In the Theorem \ref{invariants} we
summarize some of them for further usage.

\begin{theorem} \label{invariants}
Let $G$ be a finite $p$-group, and let $F$ be a field of
characteristic $p$. Then the following properties of $G$ are
determined by the group algebra $FG$:\par (i) the exponent of the
group $G$ (\cite{Kulshammer}; see also \cite{Sandling-meta});\par
(ii) the
isomorphism type of the center of the group $G$
    (\cite{Sehgal, Ward}); \par
(iii) the isomorphism type of the factorgroup $G/G'$
    (\cite{Ward}; see also \cite{Passman, Sandling-survey}) \par
(iv) the minimal number of generators $d(G')$ of the commutator
     subgroup $G'$ (follows immediately from Prop.III.1.15(ii)
     of \cite{Sehgal-topics}); \par
(v) the length of the Brauer-Jennings-Zassenhaus
$\mathcal{M}$-series of the group $G$, that is $\mathcal{M}_1(G)=G,\ \
\mathcal{M}_{n+i}(G)=(\mathcal{M}_n(G),G)\mathcal{M}_i(G)^p$,
where $i$ is the smallest integer such that $ip>n$, as well
as the isomorphism type of their factors
$\mathcal{M}_{i}(G)/\mathcal{M}_{i+1}(G)$,
$\mathcal{M}_{i}(G)/\mathcal{M}_{i+2}(G)$,
$\mathcal{M}_{i}(G)/\mathcal{M}_{2i+1}(G)$ ( \cite{PassiSehgal72, Ritter}).
\end{theorem}

It could be useful for further researches in MIP to extend the
above list of invariants determined by the modular group algebra
by adding to it at least the nilpotency class of a group. We are
able to it in several cases, listed below.

\begin{theorem} \label{Theor-class}
Let $G$ be a $p$-group and let $F$ be a field of characteristic
$p$. Then $\cl(G)$ is determined by the group algebra $FG$
provided one of the following conditions holds:\par
(i)   $\exp G=p$; \par
(ii)  $\cl(G)=2$; \par
(iii) $G'$ is cyclic; \par
(iv)  $G$ is a group of maximal class and contains an
      abelian subgroup of index $p$.
\end{theorem}

In addition, we give an application of Theorem \ref{Theor-class},
solving MIP for finite nonabelian $p$-groups containing a cyclic
subgroup of index $p^2$ (these groups were classified by Ninomia
in \cite{Ninomiya}):

\begin{theorem} \label{Theor-ninomip}
Let $G$ be a $p$-group containing a cyclic subgroup of index
$p^2$, and let $F$ be the field of $p$ elements. If for a group
$H$ we have $FG \cong FH,$ then $G \cong H$.
\end{theorem}

Our notations are standard.
$\Delta=\Delta_{K}(G)$ denotes the augmentation ideal of the
modular group algebra $KG$. $C_{2^m}$ denotes the cyclic group of order $2^m$.
We will also use the following notations for 2-groups of order $2^{m}$ and exponent $p^{m-1}$:\par
\noindent - the dihedral group $D_{m} =
  \langle a, b | a^{2^{m-1}}=1, b^2=1,           b^{-1}ab=a^{-1}         \rangle$, $m \geq 3;$\par
\noindent - the generalized quaternion group $Q_{m} =
  \langle a, b | a^{2^{m-1}}=1, b^2=a^{2^{m-2}}, b^{-1}ab=a^{-1}         \rangle$, $m \geq 3;$\par
\noindent - the semidihedral group $S_{m} =
  \langle a, b | a^{2^{m-1}}=1, b^2=1,           b^{-1}ab=a^{-1+2^{m-2}} \rangle$, $m \geq 4;$\par
\noindent - the quasi-dihedral group $M_m(2) =
  \langle a, b | a^{2^{m-1}}=1, b^2=1,           b^{-1}ab=a^{ 1+2^{m-2}} \rangle$, $m \geq 4;$\par

\section{Determination of the nilpotency class of a $p$-group}
\label{nilclass}

In this section we give the proof of the Theorem \ref{Theor-class}.

\begin{proof}
(i) Let $\mathcal{M}_i(G)$ be the $i$-th term of the
Brauer-Jennings Zassenhaus $\mathcal{M}$-series of $G$, that is
$$\mathcal{M}_1(G)=G,\ \
\mathcal{M}_{n+i}(G)=(\mathcal{M}_n(G),G)\mathcal{M}_i(G)^p,$$
where $i$ is the smallest integer such that $ip>n$. It is clear
that if $\exp G=p$, then for all positive integers $i$ we have
$\mathcal{M}_i(G)=\gamma_i(G)$. Hence, by Theorem 3(i) of
\cite{PassiSehgal72} we have
$\gamma_i(G)/\gamma_{i+1}(G)=\gamma_i(H)/\gamma_{i+1}(H)$ for all
$i\geqslant 1$. This means that $\cl(G)=cl(H)$.

(ii) It is well known that if $x$ is a non-central element of $G$
and $C_x$ is the conjugacy class of $x$ in $G$ then
$\hat{C_x}=\sum\limits_{x\in C_x}x$ lies in the subspace
$[FG,FG]$. Moreover the ideal $\langle [FG,FG]\rangle$ of $FG$ is
equal to $\Delta(G')FG$. Hence the ideal of $FG$ generated by all
central elements of $\Delta(G)$ and the subspace $[FG,FG]$ is
equal to $\Delta(Z(G)G')FG$. In particular, the order $|Z(G)G'|$
is determined by $FG$. But the orders $|Z(G)|$ and $|G'|$ are
determined, so is $|Z(G)\cap G'|=\frac{|Z(G)||G'|}{|Z(G)G'|}$.
Since $\cl(G)=2$ if and only if $G'=Z(G)\cap G'$, one can
recognize it from the structure of $FG$.

(iii) Let $G$ be a $p$-group with $G'$ cyclic of order $p^m$. Let
$H$ be a group such that $FG\cong FH$.
It follows immediately from Prop. III.1.15(ii) of \cite{Sehgal} that
$$d(G') = dim_F ( \Delta(G') FG / \Delta(G') \Delta(G) )$$
is determined by $FG$. Thus, since $G/G'$ and $d(G')$
are determined, we obtain that $H'$ is also cyclic of order $p^m$.
Let $G'= \langle g|g^{p^m}=1 \rangle $, and $H'= \langle
h|h^{p^m}=1 \rangle $. To prove the theorem, we will use induction
on the nilpotency class of $G$. For $G$ abelian the statement is
obviously true. Let $\cl(G)=c>1$ and assume that the theorem is
proved for all groups with the nilpotency class less than $c$.
Since $|Z(G) \cap G'|$ is determined, we may assume that $Z(G)
\cap G' = \langle g^{p^k} \rangle $, $Z(H) \cap H' = \langle
h^{p^k} \rangle $. Consider the ideal
$[FG,FG]FG=\Delta(G')FG=(g-1)FG$. Then
$\Delta(G')^{p^k}FG=(g^{p^k}-1)FG$, and $FG/\Delta(G')^{p^k}FG
\cong F[G/ \langle g^{p^k} \rangle ]$. Repeating the same
conclusions for $FH$, we get $F[G/ \langle g^{p^k} \rangle ] \cong
F[H/ \langle h^{p^k} \rangle ]$, or $F\overline G \cong F\overline
H$, where $\overline G = G/ \langle g^{p^k} \rangle , \overline
H=H/ \langle h^{p^k} \rangle $. Since $\cl(\overline G)=\cl(G)-1$,
$\cl(\overline H)=\cl(H)-1$, we get by induction that
$\cl(\overline G)=\cl(\overline H)$ and then $\cl(G)=\cl(H)$,
which proves the theorem.

(iv) See \cite{CJM1988}, Theorem 3.2.

\end{proof}

For $p>2$ one can give another proof of (ii). Namely, it can be
noted that in this case the ideal generated by all central
elements from the ideal $\Delta(G')FG$ is equal to $\Delta(G')FG$
if and only if $\cl(G)=2$. The example of $2$-groups of maximal
class shows that it is not the case, when $p=2$.

It is worth to note that the nilpotency class of the group of
units of $FG$ does not depend on the nilpotency class of $G$. As
it was shown in Theorem B of \cite{Shalev90}, if $G'$ is cyclic,
then the nilpotency class of $U(FG)$ is equal to $|G'|$.

\section{Presentations of finite 2-groups with a cyclic subgroup of \mbox {index 4}}
\label{presentations}

Finite nonabelian $p$-groups of order $p^m$ and exponent $p^{m-2}$
are classified in \cite{Ninomiya}. For the case $p=2$, they are
given by the following presentations: \par

\vspace{15pt} \noindent (a) $m \geq 4$: \par
$G_1 = \langle a, b | a^{2^{m-2}}=1, b^4=1, b^{-1}ab=a^{1+2^{m-3}} \rangle ; $ \par
$G_2 = Q_{m-1} \times C_2 =
       \langle a, b, c | a^{2^{m-2}}=1, b^2=a^{2^{m-3}}, c^2=1, b^{-1}ab=a^{-1},$ \par
       \hspace{250pt} $ac=ca, bc=cb \rangle$; \par
$G_3 = D_{m-1} \times C_2 =
       \langle a, b, c | a^{2^{m-2}}=1, b^2=1, c^2=1, b^{-1}ab=a^{-1},$ \par
       \hspace{250pt} $ac=ca, bc=cb \rangle$; \par
$G_4 = \langle a, b, c | a^{2^{m-2}}=1, b^2=1, c^2=1, ab=ba, ac=ca, c^{-1}bc=a^{2^{m-3}}b \rangle ; $ \par
$G_5 = \langle a, b, c | a^{2^{m-2}}=1, b^2=1, c^2=1, ab=ba, c^{-1}ac=ab, bc=cb \rangle ; $ \par

\vspace{15pt} \noindent (b) $m \geq 5$: \par
$G_6 = \langle a, b | a^{2^{m-2}}=1, b^4=1, b^{-1}ab=a^{-1} \rangle ; $ \par
$G_7 = \langle a, b | a^{2^{m-2}}=1, b^4=1, b^{-1}ab=a^{-1+2^{m-3}} \rangle ; $ \par
$G_8 = \langle a, b | a^{2^{m-2}}=1, b^4=a^{2^{m-3}}, b^{-1}ab=a^{-1} \rangle ; $ \par
$G_9 = \langle a, b | a^{2^{m-2}}=1, b^4=1, a^{-1}ba=b^{-1} \rangle ; $ \par
$G_{10} = M_{m-1}(2) \times C_2 =
\langle a, b, c | a^{2^{m-2}}=1, b^2=1, c^2=1, b^{-1}ab=a^{ 1+2^{m-3}}, $ \par
       \hspace{250pt} $ac=ca, bc=cb \rangle$; \par
$G_{11} = S_{m-1} \times C_2    =
\langle a, b, c | a^{2^{m-2}}=1, b^2=1, c^2=1, b^{-1}ab=a^{-1+2^{m-3}}, $ \par
       \hspace{250pt} $ac=ca, bc=cb \rangle$; \par
$G_{12} = \langle a, b, c | a^{2^{m-2}}=1, b^2=1, c^2=1, ab=ba, c^{-1}ac=a^{-1}, c^{-1}bc=a^{2^{m-3}}b \rangle ; $ \par
$G_{13} = \langle a, b, c | a^{2^{m-2}}=1, b^2=1, c^2=1, ab=ba, c^{-1}ac=a^{-1}b, bc=cb \rangle ; $ \par
$G_{14} = \langle a, b, c | a^{2^{m-2}}=1, b^2=1, c^2=a^{2^{m-3}}, ab=ba, c^{-1}ac=a^{-1}b, bc=cb \rangle ; $ \par
$G_{15} = \langle a, b, c | a^{2^{m-2}}=1, b^2=1, c^2=1, b^{-1}ab=a^{1+2^{m-3}}, c^{-1}ac=a^{-1+2^{m-3}},$
          \par \hspace{250pt} $bc=cb \rangle ; $ \par
$G_{16} = \langle a, b, c | a^{2^{m-2}}=1, b^2=1, c^2=1, b^{-1}ab=a^{1+2^{m-3}}, c^{-1}ac=a^{-1+2^{m-3}},$
          \par \hspace{250pt} $c^{-1}bc=a^{2^{m-3}}b \rangle ; $ \par
$G_{17} = \langle a, b, c | a^{2^{m-2}}=1, b^2=1, c^2=1, b^{-1}ab=a^{1+2^{m-3}}, c^{-1}ac=ab, bc=cb \rangle ; $ \par
$G_{18} = \langle a, b, c | a^{2^{m-2}}=1, b^2=1, c^2=b, b^{-1}ab=a^{1+2^{m-3}}, c^{-1}ac=a^{-1}b \rangle ; $ \par

\vspace{15pt} \noindent (c) $m \geq 6$: \par
$G_{19} = \langle a, b | a^{2^{m-2}}=1, b^4=1, b^{-1}ab=a^{ 1+2^{m-4}} \rangle ; $ \par
$G_{20} = \langle a, b | a^{2^{m-2}}=1, b^4=1, b^{-1}ab=a^{-1+2^{m-4}} \rangle ; $ \par
$G_{21} = \langle a, b | a^{2^{m-2}}=1, a^{2^{m-3}}=b^4, a^{-1}ba=b^{-1} \rangle ; $ \par
$G_{22} = \langle a, b, c | a^{2^{m-2}}=1, b^2=1, c^2=1, ab=ba, c^{-1}ac=a^{ 1+2^{m-4}}b$,
          \par \hspace{250pt} $c^{-1}bc=a^{2^{m-3}}b \rangle ; $ \par
$G_{23} = \langle a, b, c | a^{2^{m-2}}=1, b^2=1, c^2=1, ab=ba, c^{-1}ac=a^{-1+2^{m-4}}b$,
          \par \hspace{250pt} $c^{-1}bc=a^{2^{m-3}}b \rangle ; $ \par
$G_{24} = \langle a, b, c | a^{2^{m-2}}=1, b^2=1, c^2=1, b^{-1}ab=a^{1+2^{m-3}}, c^{-1}ac=a^{-1+2^{m-4}}b$,
          \par \hspace{250pt} $bc=cb \rangle ; $ \par
$G_{25} = \langle a, b, c | a^{2^{m-2}}=1, b^2=1, c^2=a^{2^{m-3}}, b^{-1}ab=a^{1+2^{m-3}}$,
          \par \hspace{200pt} $ c^{-1}ac=a^{-1+2^{m-4}}b, bc=cb \rangle ; $ \par
          \vspace{15pt} \noindent (d) $m = 5$: \par
$G_{26} = \langle a, b, c | a^8=1, b^2=1, c^2=a^4, b^{-1}ab=a^5, c^{-1}ac=ab , bc=cb \rangle ; $ \par

\vspace{15pt}

In the next table we listed some properties of these groups that
are important for the proof of the Theorem \ref{Theor-ninomip}.

$$
\begin{array}{|c|c|c|c|}
\hline
n   & \gamma_2(G)                        & Z(G)                                                      & \cl(G) \\
\hline
1   & C_2=\langle a^{2^{m-3}} \rangle    & C_{2^{m-3}} \times C_2 = \langle a^2, b^2 \rangle         &  2     \\
\hline
2   & C_{2^{m-3}}= \langle a^2 \rangle   & C_2 \times C_2 = \langle a^{2^{m-3}}, c \rangle           &  m-2   \\
\hline
3   & C_{2^{m-3}}= \langle a^2 \rangle   & C_2 \times C_2 = \langle a^{2^{m-3}}, c \rangle           &  m-2   \\
\hline
4   & C_2= \langle a^{2^{m-3}} \rangle   & C_{2^{m-2}} = \langle a \rangle                           &  2    \\
\hline
5   & C_2= \langle b \rangle             & C_{2^{m-3}} \times C_2 = \langle a^2,b \rangle                         &  2    \\
\hline
6   & C_{2^{m-3}}= \langle a^2 \rangle   & C_2 \times C_2 = \langle a^{2^{m-3}}, b^2 \rangle         &  m-2  \\
\hline
7   & C_{2^{m-3}}= \langle a^2 \rangle   & C_2 \times C_2 = \langle a^{2^{m-3}}, b^2 \rangle         &  m-2  \\
\hline
8   & C_{2^{m-3}}= \langle a^2 \rangle   & C_4 = \langle b^2 \rangle                                 &  m-2  \\
\hline
9   & C_2= \langle b^2 \rangle           & C_{2^{m-3}} \times C_2 = \langle a^2, b^2 \rangle &  2    \\
\hline
10  & C_2=\langle a^{2^{m-3}} \rangle    & C_{2^{m-3}} \times C_2 = \langle a^2, c \rangle &  2    \\
\hline
11  & C_{2^{m-3}}= \langle a^2 \rangle   & C_2 \times C_2 = \langle a^{2^{m-3}}, c \rangle           &  m-2  \\
\hline
%
12  & C_{2^{m-3}}= \langle a^2 \rangle   & C_4 = \langle a^{2^{m-4}}b \rangle                          &  m-2  \\
\hline
13  & C_{2^{m-3}}= \langle a^2b \rangle  & C_2 \times C_2 = \langle a^{2^{m-3}}, b \rangle           &  m-2  \\
\hline
14  & C_{2^{m-3}}= \langle a^2b \rangle  & C_2 \times C_2 = \langle a^{2^{m-3}}, b \rangle           &  m-2  \\
\hline
15  & C_{2^{m-3}}= \langle a^2b \rangle  & C_2 = \langle a^{2^{m-3}} \rangle                         &  m-2  \\
\hline
16  & C_{2^{m-3}}= \langle a^2b \rangle  & C_2 = \langle a^{2^{m-3}} \rangle                         &  m-2  \\
\hline
17  & C_2 \times C_2=\langle a^{2^{m-3}},
                              b \rangle  & C_{2^{m-4}} = \langle a^4 \rangle                 &  3    \\
\hline
18  & C_{2^{m-3}}= \langle a^2b \rangle  & C_4 = \langle a^2         \rangle, m=5                            &  m-2  \\
    &                                    & C_2 = \langle a^{2^{m-3}} \rangle, m>5                            &       \\
\hline
19  & C_4= \langle a^{2^{m-4}} \rangle   & C_{2^{m-4}} = \langle a^4 \rangle                                 &  2    \\
\hline
20  & C_{2^{m-3}}= \langle a^2 \rangle   & C_2 = \langle a^{2^{m-3}} \rangle                                 &  m-2  \\
\hline
21  & C_4=\langle b^2 \rangle            & C_{2^{m-3}}= \langle a^2 \rangle                          &  3    \\
\hline
22  & C_4= \langle a^{2^{m-4}}b \rangle  & C_{2^{m-3}} = \langle a^{2}b \rangle                                               &  3    \\
\hline
23  & C_{2^{m-3}}= \langle a^2b \rangle  & C_4 = \langle a^{2^{m-4}}b \rangle                        &  m-2  \\
\hline
24  & C_{2^{m-3}}= \langle a^2b \rangle  & C_2 = \langle a^{2^{m-3}} \rangle                         &  m-2  \\
\hline
25  & C_{2^{m-3}}= \langle a^2b \rangle  & C_2 = \langle a^{2^{m-3}} \rangle                         &  m-2  \\
\hline
26  & C_2 \times C_2=\langle a^4,
                              b \rangle  & C_2= \langle a^4 \rangle                                  &  3    \\
\hline
\end{array}
$$

\section{The Modular Isomorphism Problem for finite 2-groups\\
         containing a cyclic subgroup of index 4}
\label{mip}

In this section we give the proof of the Theorem
\ref{Theor-ninomip} for the case when $p=2$. It appears that some
of 2-groups of order $2^m$ and exponent $2^{m-2}$ are either
metacyclic groups or 2-groups of almost maximal class, for which
the modular isomorphism problem is already solved \cite{almost,
Sandling-meta}. Thus, it suffices to show that modular group
algebras of remaining 2-groups containing a cyclic subgroup of
index 4 are non-isomorphic pairwise.

\begin{proof}
We may assume that $m \geq 6$, since for $m<6$ the modular
isomorphism problem is already solved (see the review of known
results in the section \ref{intro}). Note that for $m=6,7$ it was
also solved using computer (see \cite{Bleher, Wursthorn-dip,
Wursthorn}).

Let $H$ be a finite 2-group of order $2^m$ and exponent $2^{m-2}$,
and $KG \cong KH$. Since the exponent of the group is determined
by its group algebra \cite{Sandling-meta}, it follows that $H$ is
also a 2-group of order $2^m$ and exponent $2^{m-2}$. Thus, the
family of finite 2-groups of order $2^m$ and exponent $2^{m-2}$ is
determined, and to complete the proof it remains to show that
group algebras of such groups are non-isomorphic pairwise.

First we note that $G_n$ is metacyclic for $n \in \{1, 6, 7, 8, 9,
19, 20, 21 \}$. Then these groups are determined by their modular
group algebras by \cite{Sandling-meta}.

Among the remaining non-metacyclic groups $G_n$ is a 2-group of
almost maximal class for $ n \in \{2, 3, 11, 12, 13, 14, 15, 16,
18, 23, 24, 25 \}$. As it was shown in \cite{almost}, these groups
have different sets of invariants, determined by their modular
group algebras, so their modular group algebras are non-isomorphic
pairwise.

Thus, it remains to deal with $G_n$ for $n \in \{4,5,10,17,22\}$.

Indeed, any of these groups can not be isomorphic to any of the
2-groups of almost maximal class for $ n \in \{2, 3, 11, 12, 13,
14, 15, 16, 18, 23, 24, 25 \}$. The derived subgroups of $G_{17}$
is 2-generated which splits it apart from the mentioned twelve
2-groups of almost maximal class as well as from groups $G_{4},
G_{5}, G_{10}, G_{22}$. The latter ones have the cyclic commutator
subgroup, so their nilpotency class (which is equal to 2 or 3) is
determined by the Theorem \ref{Theor-class} (iii), and this splits
them from the groups of almost maximal class as well.

Now since the isomorphism type of the center $Z(G)$ is determined
by \cite{Sehgal}, we may split the groups $G_4$ and $G_{10}$.

Groups $G_{5}$ and $G_{22}$ have isomorphic centers, but they have
cyclic commutator subgroups. Thus, again we may apply Theorem
\ref{Theor-class} (iii) and split them since $\cl(G_5)=2$ while
$\cl(G_{22})=3$, and this completes the proof for the case $p=2$.

\end{proof}

\section{The Modular Isomorphism Problem for finite $p$-groups,
$p>2$, containing a cyclic subgroup of index $p^2$} \label{mip3}

The case $p>2$ is easier and we need only the following property
of these groups:
\par

\begin{lemma}
Let $G$ be a finite nonabelian $p$-group, $p>2$, containing a
cyclic subgroup of index $p^2$. Then $G$ satisfies at least one of
the following three conditions:\par (a) $G$ is metacyclic;\par (b)
$|G'|=p$;\par (c) $G'$ is elementary abelian of order $p^2$ and
$d(G)=2$.
\end{lemma}

\begin{proof}
For $p$-groups of order $\leqslant p^4$ the lemma is clear. So
assume that $|G|=p^n$, $n>4$ and suppose $G$ is not metacyclic.
Since $\exp G=p^{n-2}$, $G$ cannot be a $p$-group of maximal class
(see for instance \cite{Leedham} 3.3). Hence by the main theorem
of \cite{Blackburn}, $G$ contains a normal subgroup $A$ of order
$p^3$ and exponent $p$. Let $C$ be a cyclic subgroup of index
$p^2$ in $G$. It is clear that $AC=G$ and $|A\cap C|=p$. Now,
$G/A$ is cyclic and $G/G'$ is not cyclic, so $|G'|$ divides $p^2$.
If $|G'|=p^2$, then obviously $2=d(G/G')=d(G)$ and the lemma
follows.
\end{proof}

Now we are able to complete the proof of the Theorem
\ref{Theor-ninomip} for the case $p>2$.

\begin{proof}

{\it Case (a).} If $G$ is metacyclic, then it is determined by
\cite{Sandling-meta}.

{\it Case (b).} If $G$ has the commutator subgroup of order $p$,
then it is determined by \cite{Sandling-elab}.

{\it Case (c).} If $d(G)=2$ and $G'$ is elementary abelian of
order $p^2$, then $G$ is determined by \cite{Baginski-elab}.

\end{proof}

\centerline{\bf Acknowledgments}
The second author wishes to thank Adalbert and
Victor Bovdi for drawing attention to the classification
\cite{Ninomiya} and their warm hospitality during his stay at the
University of Debrecen in March-April 2002, and also the NATO
Science Fellowship Programme for the support of this visit.

\bibliographystyle{amsplain}

\providecommand{\bysame}{\leavevmode\hbox to3em{\hrulefill}\thinspace}

\vspace{20pt}

\begin{minipage}[t]{7.5cm}
\small{Czes{\l}aw Bagi\'nski\\
       Institute of Computer Science\\
       Technical University of Bia{\l}ystok\\
       Wiejska 45A, 15-351\\
       Bia{\l}ystok, Poland\\
       e-mail: baginski@pb.bialystok.pl}
\end{minipage}
\begin{minipage}[t]{7.5cm}
\small{Alexander Konovalov\\
       Department of Mathematics\\
       Zaporozhye National University\\
       ul.Zhukovskogo, 66, Zaporozhye\\
       69063, Ukraine\\
       e-mail: konovalov@member.ams.org\\
       http://ukrgap.exponenta.ru/konoval.htm\\
current address: \\
Department of Mathematics\\ 
Vrije Universiteit Brussel \\
Pleinlaan 2, B-1050 Brussel, Belgium}
\end{minipage}

\end{document}